\documentclass[10pt]{article} 
\usepackage{amsfonts} 
 
\usepackage{amscd}

\newtheorem{theorem}{Theorem}[section]

\newtheorem{lemma}[theorem]{Lemma}

\newtheorem{definition}{Definition}[section]

\begin{document} 
\title{Affine structures on filiform Lie algebras} 
\author{Elisabeth Remm\\
Universit\'{e} de Haute Alsace, F.S.T.\\ 
4, rue des Fr\`{e}res Lumi\`{e}re - 68093 MULHOUSE - France} 
\date{}
\maketitle

The aim of this note is to prove that every non characteristically nilpotent
filiform algebra is provided with an affine structure. We generalize this
result to the class of nilptent algebras whose derived algebra admits non
singular derivation.

\section{Affine structures on  Lie algebras}
  
\begin{definition}
An affine connexion on a manifold $M$ is a law $\nabla $ 
which gives for every vectorfield $X$ an endomorphism $\nabla _{X}$ 
of  $\mathcal{D}^{1}(M)$ the space of vectorfiels on $M$ satisfying the
two conditions 
\begin{eqnarray*}
(1)\quad \nabla _{fX+gY} &=&f\nabla _{X}+g\nabla _{Y}; \\
(2)\quad \nabla _{X}(fY) &=&f\nabla _{X}(Y)+(Xf)Y
\end{eqnarray*}
for $f,g\in C^{\infty }(M),$ $X,Y\in \mathcal{D}^{1}(M).$\
\end{definition}

If $M$ is an Lie group $G$, then the affine connexion $\nabla $ is called
left invariant if the connexion $\nabla^{\prime}$ given by 
\[
\nabla _{X}^{\prime }(Y)=(\nabla _{X^{\Phi }}(Y^{\Phi }))^{\Phi ^{-1}}\qquad
X,Y\in \mathcal{D}^{1}(M) 
\]  
satisfies 
$$
\nabla^{\prime}=\nabla $$
for every left tanslation $\Phi$.
It is equivallent to say that the left translations are affine mapping on the
affine Lie group $(G,\nabla )$.

\begin{definition}
The torsion of the affine connexion  $\nabla$ is the tensor $T$
defined by 
$$
T(X,Y)=\nabla _{X}(Y)-\nabla _{Y}(X)-\left[ X,Y\right]
$$
The curvature of  $\nabla $ is given by the tensor $C$ defined by
$$
C(X,Y)=\nabla _{X}\nabla _{Y}-\nabla _{Y}\nabla _{X}-\nabla _{\left[
X,Y\right] } 
$$
\end{definition}
In this work we consider the left invariant connection satisfying $T=0,C=0$.

Let $\frak{g}$ be the Lie algebra of the affine group $G$. As the operator
$\nabla$ is left invariant, it inducces a mapping, always noted $\nabla$ :
$$ 
\nabla : \frak{g} \otimes \frak{g} \rightarrow  \frak{g}.
$$
If we put $X.Y= \nabla_{X}Y$ for every $X,Y \in \frak(g)$, then 
the condition on the curvature and torsion imply that this product
satisfies :  
$$
1) X.\left( Y.Z\right) -Y.\left( X.Z\right) =\left( X.Y\right)
.Z-\left(  Y.X\right) .Z 
$$ 
$$
2) X.Y-Y.X=\left[ X,Y\right]  
$$
for every $X,Y,Z\in \frak{g}$ .
\begin{definition}
A such operator $\nabla$ on the Lie algebra $\frak{g}$ is called an affine
structure on $\frak{g}$.
\end{definition}

\section{Affine structure on nilpotent Lie algebras}
The problem of existence of affine structures on nilpotent Lie algebras has
been put by John Milnor. Recently, Benoist has proposed  examples of
11-dimensional nilpotent Lie algebras which are not endowed with such
structure.  These Lie algebras have the following structure
\[
\begin{tabular}{ll}
$\left[ X_{1},X_{i}\right] =X_{i+1}\qquad i=2,..,10$ &  \\ 
$\left[ X_{2},X_{4}\right] =X_{6}$ &  \\ 
$\left[ X_{2},X_{6}\right] =-5X_{8}+2X_{9}+2tX_{10}$ &  \\ 
$\left[ X_{2},X_{8}\right] =\frac{26}{5}X_{10}+\frac{28}{25}X_{11}$ &  \\ 
$\left[ X_{3},X_{4}\right] =3X_{7}-X_{8}-tX_{9}$ &  \\ 
$\left[ X_{3},X_{6}\right] =-\frac{12}{5}X_{9}-\frac{1}{25}X_{10}+\frac{%
-448+1525t}{2000}X_{11}$ &  \\ 
$\left[ X_{3},X_{8}\right] =\frac{321}{80}X_{11}$ &  \\ 
$\left[ X_{4},X_{6}\right] =\frac{27}{5}X_{10}-\frac{24}{25}X_{11}$ &  \\ 
$\left[ X_{5},X_{6}\right] =\frac{1377}{80}X_{11}$ &  \\ 
$\left[ X_{2},X_{3}\right] =X_{5}$ &  \\ 
$\left[ X_{2},X_{5}\right] =-2X_{7}+X_{8}+tX_{9}$ &  \\ 
$\left[ X_{2},X_{7}\right] =-\frac{13}{5}X_{9}+\frac{51}{25}X_{10}+\frac{%
448+2475t}{2000}X_{11}$ &  \\ 
$\left[ X_{2},X_{9}\right] =\frac{19}{16}X_{11}$ &  \\ 
$\left[ X_{3},X_{5}\right] =3X_{8}-X_{9}-tX_{10}$ &  \\ 
$\left[ X_{3},X_{7}\right] =\frac{-39}{5}X_{10}+\frac{23}{25}X_{11}$ &  \\ 
$\left[ X_{4},X_{5}\right] =\frac{27}{5}X_{9}-\frac{24}{25}X_{10}+\frac{%
448-3525t}{2000}X_{11}$ &  \\ 
$\left[ X_{4},X_{7}\right] =-\frac{189}{16}X_{11}\qquad $ & 
\end{tabular}
t \in \mathbb{R}
\]
{\bf Classical examples}

\noindent 1.Dimensions less than 7

Every Lie nilpotent Lie algebras of dimension less or equal to $7$ admits an
affine structure. 

\noindent 2. Symplectic Lie algebras

Let $\frak{g}$ a $2p$-dimensional Lie algebra endowed to a symplectic form 
$\theta \in \Lambda ^{2}\frak{g}^{*}$. It satisfies $d\theta =0$ where 
$$
d\theta (X,Y,Z)=\theta (X,\left[ Y,Z\right] )+\theta (Y,\left[ Z,X\right]
)+\theta (Z,\left[ X,Y\right] ). 
$$
For every $X\in \frak{g}$ let $f(X)$ be defined by 
$$
\theta (adX(Y),Z)=-\theta (Y,f(X)(Z)). 
$$
Then $\nabla_{X}Y =f(X)(Y)$ is an affine structure
$\frak{g}$.

\noindent 3. Lie algebras admitting a regular derivation

Such an algebra is necessary nilpotent. More there exists a diagonalizable
regular derivation. Let $f$ be such derivation. For every 
 $X\in \frak{g}$ we put 
$$
\nabla _{X}=f^{-1}\circ adX\circ f. 
$$
Then
\begin{eqnarray*}
\nabla _{X}Y-\nabla _{Y}X &=&f^{-1}([X,f(Y)]-[Y,f(X)]) \\
&=&f^{-1}(f([X,Y])) \\
&=&[X,Y].
\end{eqnarray*}
We have also 
\begin{eqnarray*}
\nabla _{X}\nabla _{Y}-\nabla _{Y}\nabla _{X} &=&f^{-1}\circ adX\circ f\circ
f^{-1}\circ adY\circ f \\
&&-f^{-1}\circ adY\circ f\circ f^{-1}\circ adX\circ f \\
&=&f^{-1}\circ ad[X,Y]\circ f \\
&=&\nabla _{[X,Y]}.
\end{eqnarray*}
This operator $\nabla $ defines an affine structure on $\frak{g}$.

\section{Non characteristically filiform algebra}
Recall that a Lie algebra is called characteristically nilpotent is every
derivation is nilpotent. Examining the counter examples of Benoist and Burde we
the following conjecture becomes natural:

\noindent {\bf Conjecture} [Kh] Every nilpotent Lie algebra which doest not
admit affine structure is characteristically nilpotent.

Of course, the converse is false. There exist seven dimensional
characteristically nilpotent Lie algebras and these algebras have affine
structures.

The aim of this section is to prove the following result :

\begin{theorem}
Every filiform non characteristically nilpotent Lie algebra admits an affine
structure
\end{theorem}

\noindent{\bf Proof}. A $n$-dimensional nilpotent Lie algebra is called
filiform if the central descending sequence satisfies :
$$
\frak{g}\supset \mathcal{C}^{1}\frak{g}\supset \cdots \supset \mathcal{C}%
^{n-2}\frak{g}\supset \left\{ 0\right\} =\mathcal{C}^{n-1}\frak{g} 
$$
and we have
$$
\left\{ 
\begin{array}{l}
\dim \mathcal{C}^{1}\frak{g}=n-2, \\ 
\dim \mathcal{C}^{i}\frak{g}=n-i-1,\quad i=1,...,n-1.
\end{array}
\right. 
$$  
For a non characteristically nilpotent Lie algebra $\frak{g}$, let us call rank
of $\frak{g}$ the dimension of a maximal exterior torus of derivations (a
maximal abelian subalgebra of $Der(\frak{g})$ of which elements are semi-simple
derivations). We have the following results [G.K]:

1. If the Lie algebra $\frak{g}$ is filiform, its rank $%
r(\frak{g})$ satisfies 
$$
r(\frak{g})\leq 2.
$$

2.Every filiform Lie algebras of rank 2 is isomorphic to $L_{n}$ or $%
Q_{n}.$ where $L_{n}$ and $Q_{n}$ are the $n$-dimensional filiform Lie algebras
defined by  
$$
L_{n}:\left\{ [Y_{1},Y_{j}]=Y_{1+j},\quad j=2,...,n-1\right.
$$
$$
Q_{n}=\left\{ 
\begin{array}{l}
\lbrack Y_{1},Y_{j}]=Y_{1+j},\quad j=2,...,n-1 \\ 
\lbrack Y_{i},Y_{n-i+1}]=(-1)^{i+1}Y_{n},\quad i=2,...,p
\end{array}
\right. \quad n=2p.
$$
For each Lie algebra, a maximal exterior torus is precisely determined.

If $\frak{g}=L_{n}$, there exists a torus generated by the diagonal
derivations : 
$$
f_{1}(Y_{1})=0,\quad f_{1}(Y_{i})=Y_{i}, \quad 2\leq i\leq n 
$$
$$
f_{2}(Y_{1})=Y_{1},\quad f_{2}(Y_{i})=iY_{i}, \quad 2\leq i\leq n
$$
the basis $\{Y_{i}\}$ being as above.

If $\frak{g}=Q_{n}$, the basis $\{Y_{i}\}$ is not a basis of eigenvectors
for a diagonalizable derivation. We can consider the new basis given by 
$$
Z_{1}=Y_{1}-Y_{2},Z_{2}=Y_{2},...,Z_{n}=Y_{n}
$$
This basis satisfies 
$$
\lbrack Z_{1},Z_{j}] =Z_{1+j},\quad j=2,...,n-2,\quad \\
\lbrack Z_{i},Z_{n-i+1}] =(-1)^{i+1}Z_{n},\quad i=2,...,n/2
$$
Then the diagonal derivations 
$$
f_{1}(Z_{1})=0,\quad f_{1}(Z_{i})=Z_{i}, \quad 2\leq i\leq n-1,
\quad  f_{1}(Z_{n})=2Z_{n}
$$
$$
f_{2}(Z_{1})=Z_{1},\quad f_{2}(Z_{i})=(i-2)Z_{i}, \quad 2\leq i\leq
n-1,\quad f_{2}(Z_{n})=(n-3)Z_{n}.
$$
generates a maximal exterior torus of derivations.

3. Every filiform Lie algebra of rank 1 and dimension $n$ is
isomorphic to one of the following Lie algebras

i) $A_{n}^{k}\left( \lambda _{1},...,\lambda _{t-1}\right) ,$ $t=\left[ 
\frac{n-k+1}{2}\right] $ , $2\leq k\leq n-3$

$$
\left\{ 
\begin{array}{l}
\left[ Y_{1},Y_{i}\right] =Y_{i+1},\quad i=2,...,n-1 \\ 
\left[ Y_{i},Y_{i+1}\right] =\lambda _{i-1}Y_{2i+k-1}\quad , \quad 2\leq
i\leq t \\ 
\left[ Y_{i},Y_{j}\right] =a_{ij}Y_{i+j+k-2}\quad , \quad 2\leq i\leq
j\quad i+j+k-2\leq n \end{array}
\right.
$$

ii) $B_{n}^{k}\left( \lambda _{1},...,\lambda _{t-1}\right) $ $n=2m$ , $t=%
\left[ \frac{n-k}{2}\right] $ , $2\leq k\leq n-3$
$$
\left\{ 
\begin{array}{l}
\left[ Y_{1},Y_{i}\right] =Y_{i+1}\quad i=2,...,n-2 \\ 
\left[ Y_{i},Y_{n-i+1}\right] =\left( -1\right) _{n}^{i+1}Y\quad , \quad 
i=2,...,n-1 \\ 
\left[ Y_{i},Y_{i+1}\right] =\lambda _{i-1}Y_{2i+k-1}\quad ,\quad i=2,...,t
\\ 
\left[ Y_{i},Y_{j}\right] =a_{ij}Y_{i+j-k-2}\quad,\quad 2\leq i,j\leq n-2%
,i+j+k-2\leq n-2,\quad  j\neq i+1
\end{array}
\right.
$$

iii) $C_{n}^{{}}\left( \lambda _{1},...,\lambda _{t}\right) $ , $n=2m+2$ , $%
t=m-1$

$$
\left\{ 
\begin{array}{l}
\left[ Y_{1},Y_{i}\right] =Y_{i+1}\quad i=2,...,n-2 \\ 
\left[ Y_{i},Y_{n-i+1}\right] =\left( -1\right) _{n}^{i-1}Y_{n}\quad , \quad 
i=2,...m+1 \\  \left[ Y_{i},Y_{n-i-2k+1}\right] =\left( -1\right)
^{i+1}\lambda _{k}Y_{n}% \quad , \quad i=2,...,n-2-2k\quad k=1,...,m-1
\end{array}
\right.
$$

The non defined brackets are equal to zero. In this theorem, $[x]$ denotes
the integer part of $x$ and $\left( \lambda _{1},...,\lambda _{t}\right) $
are non simultaneously vanishing parameters satisfying polynomial equations
associated to the Jacobi conditions. Moreover, the constants $a_{ij}$
satisfy 
$$
a_{ij}=a_{ij+1}+a_{i+1,j}
$$
and $a_{ii+1}=\lambda _{i-1}.$

We can easily see that the filiform algebra $L_n$ , $Q_n$ or of type $A^{n}$ or
$B^{n}$ admit regular derivations. Then they admits affine structure. Let us
consider  the case $C^n$. This algebra is of rank $1$. The exterior torus of
derivation is generated by 
$$
f(Y_{1})=0,f(Y_{i})=Y_{i},\quad i=2,...,n-1,\quad f(Y_{n})=2Y_{n}. 
$$
Thus every derivation is singular. 
\begin{lemma}
The restriction of the derivation $f$ to the derived subalgebra
$D(\frak{g})$ is a regular derivation of $D(\frak{g})$.
\end{lemma}
Let us consider a vectorial endomorphism $g$ of $\frak{g}$ which leaves
invariant $D(\frak{g})$, and such that the restriction to $D(\frak{g})$
satisfies $f \circ g = Id$. Then the bilinear mapping given by 
$$
\nabla _{X}=g\circ adX\circ f. 
$$
defines an affine structure on $C^n$. In fact 
$$
\nabla _{X}(Y)-\nabla _{Y}(X) =g\circ adX\circ f(Y)-g\circ adY\circ f(X) \\
=g(f[X,Y])
$$
because $\ f$ is a derivation. As $g=f^{-1}$ on the derived subalgebra, we can
deduce  
$$
\nabla _{X}(Y)-\nabla _{Y}(X)=[X,Y]. 
$$
In the some way
$$
\nabla _{X}\nabla _{Y}(Z)-\nabla _{Y}\nabla _{X}(Z)
=g[X,[Y,f(Z)]]-g[Y,[X,f(Z)]] \\
=-g[f(Z),[X,Y]] 
$$
Then 
$$
\nabla _{X}\nabla _{Y}(Z)-\nabla _{Y}\nabla _{X}(Z)=\nabla _{[X,Y]}(Z)
$$
This proves the theorem.

\section{A theorem of existence of affine structure}
The previous proof gives the following result

\begin{theorem}
Let $\frak{g}$ be a nilpotent Lie algebra admitting a derivation of which
restriction to the derived subalgebra is regular. Then this algebra admits an
affine structure
\end{theorem}
$$
$$
\noindent BIBLIOGRAPHY

\medskip

\noindent [Be] Besnoit Y. Une nilvari\'{e}t\'{e} non affine. 
\textit{J.Diff.Geom. } \textbf{41} (1995), 21-52.

\smallskip

\noindent \noindent [Bu1] Burde D. Left invariant affine structure on
reductive Lie groups. \textit{J. Algebra}. \textbf{181}, (1996), 884-902.

\smallskip

\noindent \noindent [Bu2] Burde D. Affine structures on nilmanifolds. 
\textit{Int. J. of Math}, \textbf{7} (1996), 599-616.

\smallskip

\noindent [G.K] Goze M., Khakimdjanov Y., \textit{Nilpotent Lie algebras}.
Kluwer editor. 1995.

\smallskip

\noindent [G.R] Goze M., Remm E., Affine structures on abelian Lie algebras. Linear Algebra and its Applications. {\bf 360}, 2003. 215-230.

\smallskip

\noindent [Kh] Khakimdjanov Y.,  Characteristically nilpotent and affine Lie
algebra. \textit{Actes Colloque Vigo 2000} To appear.

\smallskip

\noindent [R] Remm E., Structures affines sur les alg\`ebres de Lie et op\'erades Lie-admissibles.
Th\`ese, Mulhouse (2001).

\smallskip

\noindent [R2] Remm E., Structures affines on contact Lie algebra. xxx: math.RA/0109077.

\end{document}